# GAUSSIAN ESTIMATES FOR SPATIALLY INHOMOGENEOUS RANDOM WALKS ON $\mathbf{Z}^D$

BY SAMI MUSTAPHA

*Institut Mathématique de Jussieu*

It is shown in this paper that the transition kernel corresponding to a spatially inhomogeneous random walk on $\mathbf{Z}^d$ admits upper and lower Gaussian estimates.

**1. Introduction.** We shall consider in this paper spatially inhomogeneous random walks $(S_j)_{j \in \mathbf{N}}$ with bounded symmetric increments in $\mathbf{Z}^d$. More precisely let $\Gamma = -\Gamma \subset \mathbf{Z}^d$ be a symmetric finite subset of $\mathbf{Z}^d$ and let $\pi : \mathbf{Z}^d \times \Gamma \longrightarrow [0,1]$ such that

$$\sum_{e \in \Gamma} \pi(x,e) = 1, \qquad \pi(x,e) = \pi(x,-e), \qquad e \in \Gamma, x \in \mathbf{Z}^d.$$

Then we let $(S_j)_{j \in \mathbf{N}}$ be the Markov chain defined by

$$\mathbf{P}[S_{j+1} = x + e // S_j = x] = \pi(x,e), \qquad e \in \Gamma, x \in \mathbf{Z}^d, j = 0, 1, \dots.$$

To avoid unnecessary complications we shall assume that $\Gamma$ contains 0 and all unit vectors in $\mathbf{Z}^d$, that is, all $e$ with $|e| = 1$ where $|\cdot|$ denotes the Euclidean norm. Furthermore we shall impose the following ellipticity condition:

$$(1.1) \qquad \pi(x,e) \geq \alpha, \qquad x \in \mathbf{Z}^d, e \in \Gamma,$$

for some $\alpha > 0$. It must be emphasized that the random walk $(S_j)_{j \in \mathbf{N}}$ is not necessarily reversible.

We shall denote by

$$(1.2) \qquad p_n(x,y) = \mathbf{P}_x[S_n = y], \qquad n = 1, 2, \dots, x, y \in \mathbf{Z}^d,$$

the transition kernel corresponding to the chain $(S_j)_{j \in \mathbf{N}}$ and by $L$ the corresponding generator, that is, the difference operator defined by

$$(1.3) \qquad Lf(x) = \sum_{e \in \Gamma} \pi(x,e)(f(x+e) - f(x)), \qquad f : \mathbf{Z}^d \to \mathbf{R}.$$

We shall prove that there exists a unique (up to a multiplicative constant) positive solution $M(\cdot)$ of the adjoint equation

$$(1.4) \quad L^*M(x) = \sum_{e \in \Gamma} \pi(x-e, e)M(x-e) - M(x) = 0, \qquad x \in \mathbf{Z}^d,$$

globally defined on $\mathbf{Z}^d$ (cf. Section 3.2). We shall denote

$$V(x, r) = \sum_{z \in B_r(x)} M(z), \qquad x \in \mathbf{Z}^d, r > 0,$$

where $B_r(x) = \{y \in \mathbf{Z}^d, |y - x| < r\}, x \in \mathbf{Z}^d, r > 0$.

THEOREM 1. *Let $(S_j)_{j \in \mathbf{N}}$ be as above. Let $p_n(x, y)$, $x, y \in \mathbf{Z}^d$, $n = 1, 2, \ldots$, denote the corresponding transition kernel. Then there exists $C > 0$, depending only on $d$, $\Gamma$ and $\alpha$, such that*

$$(1.5) \quad p_n(x, y) \leq \frac{CM(y)}{(V(x, \sqrt{n})V(y, \sqrt{n}))^{1/2}} \exp\left(-\frac{|x-y|^2}{Cn}\right),$$

$$n \geq 1, x, y \in \mathbf{Z}^d,$$

$$(1.6) \quad p_n(x, y) \geq \frac{M(y)}{C(V(x, \sqrt{n})V(y, \sqrt{n}))^{1/2}} \exp\left(-\frac{C|x-y|^2}{n}\right),$$

$$n \geq 1, x, y \in \mathbf{Z}^d, |x-y| \leq n/C.$$

The following comments may be helpful in placing the above theorem in its proper perspective.

(i) It will be shown in Section 3 that the volume function $V(x, r)$ satisfies the doubling property

$$(1.7) \qquad V(x, r) \leq CV(x, 2r), \qquad x \in \mathbf{Z}^d, r > 0.$$

The volume factor $(V(x, \sqrt{n})V(y, \sqrt{n}))^{1/2}$ in (1.5), (1.6) can therefore be replaced by $V(x, \sqrt{n})$.

(ii) In the reversible case, Theorem 1 is an immediate consequence of Delmotte's work (cf. [5]). Delmotte proved the equivalence of the upper and lower Gaussian estimates to the volume doubling property (1.7) plus the Poincaré inequality for reversible Markov chains with bounded increments on graphs. His approach relies on a clever adaptation of the Moser iteration process. The reversibility property $p(x, y)m(x) = p(y, x)m(y)$ verified by the kernel of the chain $(S_j)_{j \in \mathbf{N}}$ and its invariant measure $m$ plays a crucial role in [5] (cf. also [2]).

(iii) Reversible Markov chains on $\mathbf{Z}^d$ are the discrete analogues of diffusions generated by second-order differential operators in divergence form and the inhomogeneous walks can be considered as the analogues of second-order operators in nondivergence form (cf. [12], Table 1, page 78). The first two-sided Gaussian bound for fundamental solutions of parabolic equations in divergence form with measurable coefficients is due to Aronson (cf. [1]). For operators in nondivergence form such upper and lower Gaussian estimates were proved only recently by Escauriaza in [6].

(iv) In Aronson's work the parabolic Harnack inequality is used to obtain the Gaussian lower bound (cf. [1]). In fact both upper and lower bound for the heat kernel can easily be deduced from the parabolic Harnack principle (cf. [19]). Conversely it is shown in [10] that the two-sided Gaussian bound implies the Harnack inequality. Saloff-Coste showed in [16] that the parabolic Harnack inequality (or the two-sided Gaussian bound) for a divergence-form second-order operator (or for the Laplace–Beltrami operator on a Riemannian manifold) is equivalent to a family of Poincaré type inequalities for balls and the doubling property (cf. also [11]). The results of [5] are the discrete counterpart of [16]. It will be interesting to discuss this type of equivalence for both nondivergence differential operators and nonreversible random walks.

A general outline of the paper is as follows. Section 2 collects the main potential theoretic properties of spatially inhomogeneous random walks on $\mathbf{Z}^d$. The two new results of this section (Theorems 4 and 5) are of independent interest and are proved in Section 5. In Section 3 we define the concept of normalized adjoint solution adapted to spatially inhomogeneous random walks and we prove that adjoint solutions verify a parabolic Harnack principle. This Harnack principle is used in Section 4 to deduce the Gaussian estimates of Theorem 1.

**2. Potential theory.** Let $A \subset \mathbf{Z}^d$ denote a bounded domain (i.e., a finite connected set of vertices in $\mathbf{Z}^d$). We let

$$\partial A = \{x \in A^c, x = z + e, \text{ for some } z \in A \text{ and } e \in \Gamma\},$$

$\Gamma$ being as in the previous section, and

$$\overline{A} = A \cup \partial A.$$

Let $B = A \times \{a \leq k \leq b\} \subset \mathbf{Z}^d \times \mathbf{Z}$ where $A \subset \mathbf{Z}^d$ and where $a < b \in \mathbf{Z}$. We let

$$\partial_l B = \bigcup_{a < k < b} \partial A \times \{k\},$$
$$\partial_p B = \partial_l B \cup (\overline{A} \times \{a\}),$$

and $\overline{B} = B \cup \partial_p B$. $\partial_p B$ is the parabolic boundary of $B$ and $\partial_l B$ is its lateral boundary. We say that $u : \overline{A} \longrightarrow \mathbf{R}$ is harmonic in $A \subset \mathbf{Z}^d$ if

$$Lu(x) = \sum_{e \in \Gamma} \pi(x, e)(u(x+e) - u(x)) = 0, \qquad x \in A.$$

Let $B = A \times \{a \leq k \leq b\} \subset \mathbf{Z}^d \times \mathbf{Z}$ and $u : \overline{B} \longrightarrow \mathbf{R}$. We say that $u$ is caloric in $B$ if

$$\mathcal{L}u(x, k) = \sum_{e \in \Gamma} \pi(x, e) u(x+e, k) - u(x, k+1) = 0,$$

$$(x, k) \in A \times \{a \leq k < b\}.$$

The following maximum principle is immediate.

THEOREM 2 (Maximum principle). *Let $B = A \times \{a \leq k \leq b\} \subset \mathbf{Z}^d \times \mathbf{Z}$, where $a < b \in \mathbf{Z}$ and $A$ is a bounded domain in $\mathbf{Z}^d$ and let $u : \overline{B} \longrightarrow \mathbf{R}$ such that $\mathcal{L}u = 0$ in $B$ and $u \geq 0$ on $\partial_p B$. Then $u \geq 0$ in $B$.*

The following theorem (cf. [13]) is a random walk version of a well-known and fundamental result in the potential theory of second-order equations in nondivergence form (for an elliptic version cf. [14]).

THEOREM 3 (Parabolic Harnack principle). *Let $u$ be a nonnegative caloric function in $B_{2r}(y) \times \{s - 4r^2 \leq k \leq s\}$, $(y, s) \in \mathbf{Z}^d \times \mathbf{Z}$, $r \geq 1$. Then*

$$\sup\{u(x, k); x \in B_r(y), s - 3r^2 < k < s - 2r^2\}$$
(2.1)
$$\leq C \inf\{u(x, k); x \in B_r(y), s - r^2 < k < s\},$$

*where $C = C(d, \alpha, \Gamma) > 0$.*

In the proof of Theorem 1, together with the previous results we need the following estimates which describe the boundary behavior of nonnegative caloric functions. Let $y_0 \in \mathbf{Z}^d$, let $R_0 > 0$ and let $\Omega = B_{R_0}(y_0)$.

Let $Q = \Omega \times \mathbf{Z}$, let $Y = (y, s) \in \partial \Omega \times \mathbf{Z}$ and let $c \leq r \leq R_0/2$ where $c > 0$ denotes a constant times $\text{diam}(\Gamma)$. We shall denote

$$C_r(Y) = B_r(y) \times \{s - r^2 \leq k \leq s\}, \qquad Q_r(Y) = Q \cap C_r(Y),$$
$$\overline{Y}_r = (y_r, s + 2[r]^2), \qquad \underline{Y}_r = (y_r, s - 2[r]^2),$$

where $y_r \in \Omega$ satisfies $|y_r - (R - r/2)\frac{y - y_0}{|y - y_0|}| \leq 1$ and $[r]$ denotes the greatest integer $\leq r$.

THEOREM 4 (Boundary Harnack principle). *Let $Y = (y,s) \in \partial\Omega \times \mathbf{Z}$. Let $c \leq r \leq R_0/K$ where $K > 0$ is large enough. Assume that $u$ and $v$ are two nonnegative caloric functions in $Q \cap (B_{3Kr}(y) \times \{s - 9K^2r^2 \leq k \leq s + 9K^2r^2\})$ and $u = 0$ on $(\partial\Omega \times \mathbf{Z}) \cap (B_{2Kr}(y) \times \{s - 4K^2r^2 \leq k \leq s + 4K^2r^2\})$. Then*

$$(2.2) \qquad \sup_{Q_r(Y)} \frac{u}{v} \leq C \frac{u(\overline{Y}_{Kr})}{v(\underline{Y}_{Kr})},$$

*where $C = C(d, \alpha, \Gamma) > 0$.*

THEOREM 5 (Backward Harnack principle). *Let $u$ be a nonnegative caloric function in $B_r(y_0) \times \mathbf{N}$ ($y_0 \in \mathbf{Z}^d$, $r > 0$ large enough) vanishing on $\partial B_r(y_0) \times \mathbf{N}$. Then*

$$(2.3) \quad u(x, k + 2[r]^2) \leq Cu(x,k), \qquad (x,k) \in B_r(y_0) \times \{r^2 \leq k \leq 3r^2\},$$

*where $C = C(d, \alpha, \Gamma) > 0$.*

The proofs of the estimates (2.2) and (2.3) which follow strongly the proofs of the corresponding facts about the boundary behavior of nonnegative solutions of second-order equations in nondivergence form (cf. [3, 4, 7, 8, 9, 15]) are given in Section 5.

We shall use throughout the usual convention $f \approx g$ to indicate that $C^{-1} \leq f/g \leq C$ for an appropriate constant $C > 0$ and $C$, $c$ are used to denote different positive constants which depend only on $d$, $\alpha$, diam($\Gamma$).

**3. The adjoint Harnack principle.** Let $D \subset \mathbf{Z}^d$ denote a bounded domain and $a < b \in \mathbf{Z}$. We say that $v = v(x,t) : \overline{D \times \{a \leq t \leq b\}} \to \mathbf{R}$ is a parabolic adjoint solution of $L$ in $D \times \{a \leq t \leq b\}$, if $v$ satisfies the equation

$$v(y, t+1) - v(y,t) = L^*v(y,t), \qquad t = a, \ldots, b-1, y \in D,$$

where $L^*$ is defined as in (1.4). Let $m(\cdot)$ be a fixed positive adjoint solution for $L$ in $D$ [i.e., $m : \overline{D} \to \mathbf{R}$, $m(x) > 0$, $\forall x \in D$, $L^*m = 0$ in $D$]. For instance, if $D \subset B_{r_0}(0)$ lies in the Euclidean ball centered at 0 and of radius $r_0 > 0$ large enough, we can set $m(x) = G(x^*, x)$ where $x^* \in B_{4r_0}(0) \setminus B_{3r_0}(0)$, with $G(\cdot, \cdot)$ being the Green function of $(S_j)_{j \in \mathbf{N}}$ in the ball $B_{5r_0}(0)$. Let $v$ be a parabolic adjoint solution in $D \times \{a \leq t \leq b\}$; the function

$$\tilde{v}(y,t) = \frac{v(y,t)}{m(y)}, \qquad (y,t) \in \overline{D \times \{a \leq t \leq b\}},$$

is called a normalized parabolic adjoint solution of $L$ in $D \times \{a \leq t \leq b\}$ (cf. [4]).

THEOREM 6. *Suppose that $\tilde{v}$ is a nonnegative normalized adjoint solution for $L$ in $B_r(y_0) \times \mathbf{N}$, where $y_0 \in \mathbf{Z}^d$ and $r > 0$ large enough. Then there exists a constant $C > 0$ depending only on $d$, $\alpha$, $\Gamma$ such that*

$$\sup\{\tilde{v}(y,s); y \in B_{r/2}(y_0), r^2 < s < 2r^2\}$$
(3.1)
$$\leq C \inf\{\tilde{v}(y,s); y \in B_{r/2}(y_0), 3r^2 < s < 4r^2\}.$$

PROOF. Let $B_r(y_0) = B_r$ and let $q_t(\cdot,\cdot)$ denote the Green function of $\mathcal{L}$ in $B_r$. An easy induction on $t$ gives the following representation formula for the parabolic normalized adjoint solutions:

$$\tilde{v}(y,t) = \sum_{x \in B_r} m(x)\tilde{v}(x,0)\frac{q_t(x,y)}{m(y)}$$
(3.2)
$$+ \sum_{s=0}^{t-1} \sum_{x \in \partial B_r} m(x)\tilde{v}(x,s) \sum_{e \in \Gamma_x} \pi(x,e)\frac{q_{t-s-1}(x+e,y)}{m(y)},$$
$$y \in B_{r/2}, t = 1, 2, \ldots,$$

where $\Gamma_x = \{e \in \Gamma / x + e \in B_r\}$. On the other hand, let us observe that if we extend $q_s(\cdot,y)$ by $q_s(\cdot,y) \equiv 0$, $s \leq 0$, in an appropriate neighborhood of $\partial B_r$ and use the boundary Harnack principle (2.2) to compare $u_1(x,t) = q_t(x,y)$ and $u_2(x,t) = q_{t+10r^2}(x,y)$, combined with the backward Harnack principle (2.3), we deduce then that

$$q_s(x,y) \leq Cq_{s+r^2}(x,y), \qquad \operatorname{dist}(x, \partial B_r) \leq cr,$$
(3.3)
$$y \in B_{r/2}, 0 < s < 2r^2.$$

Let now $y_1, y_2 \in B_{r/2}$ and $r^2 < t_1 < 2r^2$, $3r^2 < t_2 < 4r^2$. By (3.2), (3.3) and Theorem 5 we have

$$\tilde{v}(y_1, t_1) \leq \sum_{x \in B_r} m(x)\tilde{v}(x,0)\frac{q_{t_1}(x,y_1)}{m(y_1)}$$
(3.4)
$$+ C\sum_{s=0}^{2r^2} \sum_{x \in \partial B_r} m(x)\tilde{v}(x,s) \sum_{e \in \Gamma_x} \pi(x,e)\frac{q_{r^2}(x+e,y_1)}{m(y_1)},$$

$$\tilde{v}(y_2, t_2) \geq \sum_{x \in B_r} m(x)\tilde{v}(x,0)\frac{q_{t_2}(x,y_2)}{m(y_2)}$$
(3.5)
$$+ \frac{1}{C}\sum_{s=0}^{2r^2} \sum_{x \in \partial B_r} m(x)\tilde{v}(x,s) \sum_{e \in \Gamma_x} \pi(x,e)\frac{q_{r^2}(x+e,y_2)}{m(y_2)}.$$

A simple use of the boundary Harnack principle [combined with (2.1) and (2.3)] allows us to deduce then

$$\tilde{v}(y_1, t_1) \leq C \frac{q_{r^2}(y_0, y_1)}{m(y_1)} \frac{m(y_2)}{q_{r^2}(y_0, y_2)} \tilde{v}(y_2, t_2).$$

Thus, to prove (3.1) it suffices to show that

(3.6) $$\frac{q_{r^2}(y_0, y_1)}{m(y_1)} \approx \frac{q_{r^2}(y_0, y_2)}{m(y_2)}, \qquad y_1, y_2 \in B_{r/2}.$$

The estimate (3.6) is a consequence of (3.4), (3.5) and the fact that 1 is a normalized parabolic adjoint solution. Indeed, let $\psi(x) = G(x, y_0)$ be the Green function of $L$ in $B_r$ with pole at $y_0$. Let $A \in B_r$ fixed with $|A - y_0| \approx r/4$. By the boundary Harnack principle and the backward Harnack principle we have

$$q_{r^2}(x, y) \approx \psi(x) \frac{q_{r^2}(y_0, y)}{\psi(A)}, \qquad x \in B_r \setminus B_{3r/4}, y \in B_{r/2}.$$

If we apply now (3.4) and (3.5) to $\tilde{v} \equiv 1$, we deduce then

$$m(y_1) \leq C \sum_{x \in B_r - B_{3r/4}} m(x) \frac{\psi(x)}{\psi(A)} q_{r^2}(y_0, y_1) + C \sum_{x \in B_{3r/4}} m(x) q_{r^2}(y_0, y_1)$$
$$+ C r^2 \sum_{x \in \partial B_r} m(x) \sum_{e \in \Gamma_x} \pi(x, e) q_{r^2}(y_0, y_1) \frac{\psi(x+e)}{\psi(A)},$$

$$m(y_2) \geq \frac{1}{C} \sum_{x \in B_r - B_{3r/4}} m(x) \frac{\psi(x)}{\psi(A)} q_{r^2}(y_0, y_2) + \frac{1}{C} \sum_{x \in B_{3r/4}} m(x) q_{r^2}(y_0, y_2)$$
$$+ \frac{r^2}{C} \sum_{x \in \partial B_r} m(x) \sum_{e \in \Gamma_x} \pi(x, e) q_{r^2}(y_0, y_2) \frac{\psi(x+e)}{\psi(A)},$$

and these two inequalities imply (3.6). □

3.1. *The doubling property for the adjoint solutions.* We start with a doubling property for the Green functions.

PROPOSITION 1. *Let $R$ be large enough and let $x_0 \in \mathbf{Z}^d$. Let $G^R(\cdot, \cdot)$ denote the Green function of $(S_j)_{j \in \mathbf{N}}$ in the ball $B_{4R}(x_0)$. There exists then a constant $C > 0$ (independent of $x_0$ and $R$) such that*

(3.7) $$\sum_{y \in B_{2r}(x_0)} G^R(x, y) \leq C \sum_{y \in B_r(x_0)} G^R(x, y), \qquad x \in B_{4R}(x_0), 1 \leq r \leq R/2.$$

This proposition is a consequence of the parabolic Harnack principle and the following lemma.

LEMMA 1. *Let $R$ be large enough and let $x_0 \in \mathbf{Z}^d$. Let $h_t^R(x,y)$, $t = 0, 1, \ldots$, $x, y \in B_{4R}(x_0)$, denote the heat kernel of $(S_j)_{j \in \mathbf{N}}$ in the ball $B_{4R}(x_0)$. Then there exists $c > 0$ (independent of $x_0$ and $R$) such that*

$$(3.8) \quad \inf_{z \in B_r(x_0)} \sum_{y \in B_{2r}(x_0)} h_s^R(z, y) \geq c, \quad 1 \leq s \leq r^2, 1 \leq r \leq R/2.$$

Indeed, by the parabolic Harnack principle we have

$$C \inf_{z \in B_{2r}(x_0)} \sum_{y \in B_r(x_0)} h_{2r^2}^R(z, y) \geq \sup_{z \in B_{2r}(x_0)} \sum_{y \in B_r(x_0)} h_{r^2}^R(z, y)$$

$$\geq \inf_{z \in B_{r/2}(x_0)} \sum_{y \in B_r(x_0)} h_{r^2}^R(z, y) \geq c.$$

We have then

$$\sum_{y \in B_r(x_0)} h_{t+2r^2}^R(x, y) \geq \sum_{y \in B_r(x_0)} \sum_{u \in B_{2r}(x_0)} h_t^R(x, u) h_{2r^2}^R(u, y)$$

$$\geq \sum_{u \in B_{2r}(x_0)} \inf_{u \in B_{2r}(x_0)} \left( \sum_{y \in B_r(x_0)} h_{2r^2}^R(u, y) \right) h_t^R(x, u)$$

$$\geq c \sum_{u \in B_{2r}(x_0)} h_t^R(x, u)$$

and this implies the doubling property for the heat kernel with the time shifted $t \to t + 2r^2$. If we sum on $t$ we deduce that

$$\sum_{y \in B_{2r}(x_0)} G^R(x, y) \leq C \sum_{y \in B_r(x_0)} \sum_{t=2r^2}^{\infty} h_t^R(x, y)$$

$$\leq C \sum_{y \in B_r(x_0)} G^R(x, y).$$

The adjoint Harnack principle (3.1) and Proposition 1 give

THEOREM 7. *Let $L^* m = 0$, $m \geq 0$, in $B_{4r}(z)$, $z \in \mathbf{Z}^d$ and $r > 0$ large enough. Then*

$$(3.9) \quad \sum_{y \in B_{2r}(z)} m(y) \leq C \sum_{y \in B_r(z)} m(y),$$

*where $C = C(d, \alpha, \Gamma)$.*

Indeed let $x^* \in B_{6r}(z) \setminus B_{5r}(z)$. Let $G(\cdot,\cdot)$ denote the Green function of $(S_j)_{j \in \mathbf{N}}$ in the ball $B_{7r}(z)$. We have

$$\sum_{y \in B_{2r}(z)} m(y) \leq C \sup_{y \in B_{2r}(z)} \frac{m(y)}{G(x^*, y)} \sum_{y' \in B_{2r}(z)} G(x^*, y')$$

$$\leq C \inf_{y \in B_{2r}(z)} \frac{m(y)}{G(x^*, y)} \sum_{y' \in B_{2r}(z)} G(x^*, y')$$

$$\leq C \sum_{y' \in B_r(z)} \frac{m(y')}{G(x^*, y')} G(x^*, y').$$

The second inequality follows from the adjoint Harnack principle and the third one from the doubling property (3.7).

PROOF OF LEMMA 1. Let

$$U(x, s) = \sum_{y \in B_{2r}(x_0)} h_s^R(x, y), \qquad x \in B_{4r}(x_0), s \geq 0.$$

Let $V(x, s)$ be the caloric function defined by

$$V(x, s+1) = \sum_{e \in \Gamma} \pi(x, e) V(x+e, s) \qquad \text{in } B_{2r}(x_0) \times \{-4r^2 \leq k < 4r^2\},$$

$$V(x, s) = 1 \qquad \text{on } \partial_p(B_{2r}(x_0) \times \{-4r^2 \leq k < 4r^2\}) \cap \{s \leq 0\},$$

$$V(x, s) = 0 \qquad \text{on } \partial_p(B_{2r}(x_0) \times \{-4r^2 \leq k < 4r^2\}) \cap \{s \geq 1\}.$$

By the maximum principle we get

$$U(x, s) \geq V(x, s), \qquad s \geq 0, x \in B_{2r}(x_0).$$

Using the parabolic Harnack inequality applied to $V(x, s)$ in $B_{2r}(x_0) \times \{-4r^2 \leq k < 4r^2\}$ we deduce that

$$\inf_{z \in B_r(x_0)} U(z, s) \geq cV(x_0, -r^2) = c, \qquad 0 < s < r^2.$$

[Note that $V \equiv 1$ on $B_{2r}(x_0) \times \{-4r^2 \leq k < 0\}$.] □

3.2. *The global adjoint solution.*

THEOREM 8. *There exists a positive adjoint solution $M$ defined globally in $\mathbf{Z}^d$. This solution is unique up to a multiplicative constant and verifies*

(3.10) $$\sum_{y \in B_{2r}(x)} M(y) \leq C \sum_{y \in B_r(x)} M(y), \qquad x \in \mathbf{Z}^d, r > 0,$$

*where $C = C(d, \alpha, \Gamma)$.*

PROOF. Let
$$m_l(y) = \alpha_l[G_{l+1}(0,y) - G_l(0,y)], \qquad y \in B_{2^l}(0), l = 1, 2, \ldots,$$
where $G_l(0, \cdot)$ is the Green function of $(S_j)_{j \in \mathbf{N}}$ in the ball $B_{2^l}(0)$, $l = 1, 2, \ldots$, with pole at the origin and where the $\alpha_l$ are chosen so that

(3.11) $$m_l(0) = 1, \qquad l = 1, 2, \ldots.$$

It is easy to see that the ellipticity condition (1.1) implies a local Harnack principle for the nonnegative adjoint solutions. This local Harnack principle and the normalization condition (3.11) imply that the $m_l$ verify
$$m_l(y) \leq C, \qquad y \in B_{2^k}(0), l \geq k,$$
with a constant $C = C(k)$ depending only on $k$. The diagonal process allows us then to deduce the existence of a global positive adjoint solution $M$ defined on $\mathbf{Z}^d$. The fact that this global adjoint solution is unique (up to a multiplicative constant) follows from the normalized adjoint Harnack principle (3.1) applied to $M_1/M_2$ where $M_1$ and $M_2$ denote two global positive adjoint solutions. Indeed it is always possible to suppose that $\inf_{\mathbf{Z}^d} M_1/M_2 = 0$ and associate to $\varepsilon > 0$, $z_\varepsilon \in \mathbf{Z}^d$, such that $M_1(z_\varepsilon)/M_2(z_\varepsilon) < \varepsilon$. By (3.1) $\sup_{B_R(z_\varepsilon)} M_1/M_2 < C\varepsilon$ with a constant $C > 0$ independent of $R$ and it suffices to let $R \to \infty$ and $\varepsilon \to 0$ to deduce that this function is constant. Finally, the doubling property (3.10) follows from Theorem 7. This completes the proof of Theorem 8. $\square$

**4. The Gaussian estimates.** The first step in proving the upper Gaussian estimate (1.5) is to prove the following mass escape estimate for $(S_j)_{j \in \mathbf{N}}$ (cf. [17, 18]).

LEMMA 2. *Let $(S_n)_{n \in \mathbf{N}}$ be as in Section 1. Let $p_n(x, y)$, $x, y \in \mathbf{Z}^d$, $n = 1, 2, \ldots$, denote the corresponding transition kernel. Then there exist $C$, $c > 0$, such that*

(4.1) $$\sum_{|x-y|>R} p_n(x,y) \leq C \exp\left(-c\frac{R^2}{n}\right), \qquad x \in \mathbf{Z}^d, n, R = 1, 2, \ldots.$$

PROOF. To prove (4.1) it suffices to prove that

(4.2) $$\sum_{l(x-y)>R} p_n(x,y) \leq e^{-cR^2/n}$$

for every linear form $l: \mathbf{R}^d \to \mathbf{R}$ such that $l(x-y) \leq |x-y|, x, y \in \mathbf{R}^d$. Let $s > 0$; we have
$$\sum_{l(x-y)>R} p_n(x,y) \leq \sum_{l(x-y)>R} e^{-sR+sl(x-y)} p_n(x,y)$$
$$\leq e^{-sR} \sum_{y \in \mathbf{Z}^d} e^{sl(x)} p_n(x,y) e^{-sl(y)}$$

from which it follows that

$$\sum_{l(x-y)>R} p_n(x,y)$$

$$\leq e^{-sR} \sum_{y_1,\ldots,y_{n-1}\in\mathbf{Z}^d} (e^{sl(x)}p_1(x,y_1)e^{-sl(y_1)})$$

(4.3)
$$\times (e^{sl(y_1)}p_1(y_1,y_2)e^{-sl(y_2)}) \times \cdots$$
$$\times (e^{sl(y_{n-1})}p_1(y_{n-1},y)e^{-sl(y)}).$$

On the other hand, we have

$$\sum_{y'\in\mathbf{Z}^d} e^{sl(y)}p_1(y,y')e^{-sl(y')} = \sum_{e\in\Gamma} \pi(y,e)e^{-sl(e)}$$

$$= 1 - s\sum_{e\in\Gamma} \pi(y,e)l(e) + O(s^2 e^{Cs})$$

$$= 1 + O(s^2 e^{Cs}),$$

where the last equality follows from the fact that $(S_n)_{n\in\mathbf{N}}$ has symmetric increments. We deduce then that

(4.4) $$\sup_{y\in\mathbf{Z}^d}\left|\sum_{y'\in\mathbf{Z}^d} e^{sl(y)}p_1(y,y')e^{-sl(y')}\right| \leq 1 + Cs^2 \leq e^{Cs^2}, \qquad |s|\leq 1.$$

If $|s|\geq 1$, it suffices to observe that

$$\left|\sum_{y'\in\mathbf{Z}^d} e^{sl(y)}p_1(y,y')e^{-sl(y')}\right| \leq \sum_{e\in\Gamma} \pi(y,e)e^{s|e|}$$

(4.5)
$$\leq e^{Cs} \leq e^{Cs^2}, \qquad y\in\mathbf{Z}^d.$$

Putting together (4.3), (4.4) and (4.5) we deduce that

$$\sum_{l(x-y)>R} p_n(x,y) \leq e^{-sR+Cns^2},$$

and optimizing over $s$, we deduce (4.2). $\square$

The second step is to apply the parabolic adjoint Harnack principle to

$$\tilde{v}(y,t) = \frac{p_t(x,y)}{M(y)}, \qquad y\in\mathbf{Z}^d, t=1,2,\ldots.$$

This gives

$$\tilde{v}(y,n) \leq C \inf_{z\in B_{\sqrt{n}}(y)} \tilde{v}(z,2n), \qquad n\geq C.$$

Hence (with the notation of Section 1)

$$\tilde{v}(y,n)V(y,\sqrt{n}) \leq C \sum_{z \in B_{\sqrt{n}}(y)} \tilde{v}(z,2n)M(z) \leq \sum_{z \notin B_{c|x-y|}(x)} p_{2n}(x,z), \qquad n \geq C,$$

in the case $\sqrt{n} < c|x-y|$, with $c > 0$ small enough. This implies, in this case, by Lemma 2

$$\tilde{v}(y,n)V(y,\sqrt{n}) \leq C \exp\left(-c\frac{|x-y|^2}{n}\right), \qquad n \geq C.$$

In the case $\sqrt{n} > c|x-y|$ it suffices to observe that the Gaussian factor in (1.5) is $\approx 1$. Hence

$$p_n(x,y) \leq \frac{CM(y)}{V(y,\sqrt{n})} \exp\left(-c\frac{|x-y|^2}{n}\right), \qquad n \geq C.$$

The doubling property (3.10) allows us to symmetrize the volume factor in this estimate and obtain

$$p_n(x,y) \leq \frac{CM(y)}{(V(x,\sqrt{n})V(y,\sqrt{n}))^{1/2}} \exp\left(-c\frac{|x-y|^2}{n}\right), \qquad n \geq C.$$

Let us observe that for $1 \leq n \leq C$ (1.5) and (1.6) are immediate consequences of the local Harnack estimate. This completes the proof of the upper estimate in Theorem 1. Finally the lower Gaussian estimate (1.6) can be deduced from the upper Gaussian estimate (1.5) and the parabolic adjoint Harnack by a standard procedure. We first use (1.5) to deduce that for $A > 0$ large enough

$$(4.6) \qquad \sum_{|x-y| \leq A\sqrt{n}} p_n(x,y) \geq \tfrac{1}{2}.$$

Parabolic adjoint Harnack applied to

$$\tilde{u}(y,t) = \frac{p_t(x,y)}{M(y)}, \qquad y \in \mathbf{Z}^d, t = 1, 2, \ldots,$$

estimate (4.6) and the doubling property (3.10) imply therefore that

$$(4.7) \qquad p_n(x,x) \geq \frac{M(x)}{CV(x,\sqrt{n})}, \qquad x \in \mathbf{Z}^d, n \geq C.$$

The lower off-diagonal estimate (1.6) is easily deduced from (4.7) by applying successively the parabolic adjoint Harnack inequality. More precisely, let us fix $x$ and $n$ as in (4.7) and let $y \in \mathbf{Z}^d$ such that $|y - x| \leq n/C$ with $C > 0$ large enough. Let $k$ be the smallest integer $\geq |x-y|^2/n$. Put

$$(a_j, t_j) = \left(a_j, \left(1 + \frac{j}{k}\right)n\right), \qquad j = 0, \ldots, k,$$

with
$$a_0 = x, \qquad a_k = y, \qquad |a_{j+1} - a_j| \approx \frac{|x-y|}{k}, \qquad 0 \leq j \leq k-1.$$

Then $(x, n) = (a_0, n)$ and $(y, 2n) = (a_k, 2n)$. Moreover
$$|a_{j+1} - a_j|^2 \approx \frac{n}{k}, \qquad 0 \leq j \leq k-1.$$

Hence the parabolic adjoint Harnack inequality yields
$$\frac{p_{2n}(x,y)}{M(y)} \geq C^k \frac{p_n(x,x)}{M(x)} \geq \frac{1}{CV(x,\sqrt{n})} \exp\left(-c\frac{|x-y|^2}{n}\right), \qquad n \geq C.$$

This completes the proof of Theorem 1.

**5. Proofs of results.** To the process $(S_j)_{j \in \mathbf{N}}$ we shall associate the corresponding space–time process
$$\dot{S}_j = (S_j, t_0 - j) \in \mathbf{Z}^d \times \mathbf{Z}, \qquad t_0 \in \mathbf{Z}, \ldots, j = 0, 1, \ldots.$$

For any cylinder $Q = \Omega \times \{a \leq k \leq b\}$, $\Omega$ being a bounded domain in $\mathbf{Z}^d$, we shall denote $\dot{\tau}_Q$ the first exit time of $\dot{S}_j$ from $Q$. The caloric measure in $Q$ at $(x_0, t_0)$ is defined by
$$\omega_Q^{(x_0,t_0)}(E) = \mathbf{P}_{(x_0,t_0)}[\dot{S}_{\dot{\tau}_Q} \in E], \qquad E \subset \partial_p Q.$$

Observe that for each $\varphi : \partial_p(\Omega \times \{a \leq k \leq b\}) \longrightarrow \mathbf{R}$, the solution of the boundary value problem
$$u(x, t+1) = \sum_{e \in \Gamma} \pi(x, e) u(x+e, t) \qquad \text{in } \Omega \times \{a \leq k < b\},$$
$$u(x,t) = \varphi(x,t) \qquad \text{on } \partial_p(\Omega \times \{a \leq k \leq b\}),$$

can be represented by means of $\omega^{x_0,t_0} = \omega_Q^{(x_0,t_0)}$, $(x_0, t_0) \in \Omega \times \{a < k \leq b\}$ as follows:
$$u(x_0, t_0) = E_{(x_0,t_0)}[\varphi(\dot{S}_{\dot{\tau}_Q})] = \sum_{(y,s) \in \partial_p Q} \varphi(y,s) \omega^{x_0,t_0}(y,s).$$

5.1. *A lower estimate for the caloric measure.* Let the notation be as in Section 2 and let $Q = \Omega \times \mathbf{Z}$. For $Y = (y, s) \in \partial_l Q$, $r > 0$ we shall denote
$$\Delta_r(Y) = \partial_l Q \cap C_r(Y),$$
where $\partial_l Q = \partial \Omega \times \mathbf{Z}$.

LEMMA 3. *Let* $\omega^X = \omega_{Q_{2r}(Y)}^{(x,t)}$. *Then*

(5.1) $$\inf_{X \in Q_r(Y)} \omega^X(\Delta_{2r}(Y)) \geq \theta, \qquad Y \in \partial_l Q, c \leq r \leq R_0/2,$$

*where* $\theta = \theta(d, \alpha, \Gamma) > 0$.

PROOF. Let $Y = (y, s) \in \partial_l Q$. It is clear that there exists a cylinder $\mathcal{C}' = B_{\mu r}(z) \times \{s - 4r^2 \leq k \leq s\} \subset C_{2r}(Y) \setminus (\Omega \times \mathbf{Z})$ (provided that $\mu$ is small enough). Let $\Delta' = B_{\mu r}(z) \times \{s - 4r^2\}$ denote the bottom of this cylinder. Using the maximum principle we deduce

$$\omega^X(\Delta_{2r}(Y)) \geq v(X) = \omega^X_{C_{2r}(Y)}(\Delta') \tag{5.2}$$

in $Q_{2r}(Y)$, and

$$v(X) \geq v'(X) = \omega^X_{\mathcal{C}'}(\Delta')$$

in $\mathcal{C}'$. On the other hand, the parabolic Harnack principle applied to $v$ gives

$$\begin{aligned}
\inf_{X \in Q_r(Y)} \omega^X(\Delta_{2r}(Y)) &\geq \inf_{X \in Q_r(Y)} v(X) \\
&\geq cv(z, s - 2r^2) \\
&\geq c'v'(z, s - 2r^2).
\end{aligned} \tag{5.3}$$

But $v'$ can be extended from $\mathcal{C}'$ to a large cylinder $\mathcal{C}'' = B_{\mu r}(z) \times [s - 6r^2, s]$ by

$$v'(X) = \omega^X_{\mathcal{C}''}(\partial_p \mathcal{C}'' \cap \{t \leq s - 4r^2\}) \tag{5.4}$$

so that $v' \equiv 1$ on $\mathcal{C}'' \cap \{t \leq s - 4r^2\}$ and the lower estimate (5.1) follows then from the Harnack principle. $\square$

COROLLARY 1. *Under the assumptions of Lemma 3, let $u$ be a non-negative solution of $\mathcal{L}u = 0$ in $Q_{3r}(Y)$, which vanishes on $\Delta_{2r}(Y)$. Then $M_r = \sup_{Q_r(Y)} u$ satisfies*

$$M_r \leq \rho M_{2r}, \qquad c \leq r \leq R_0/4, \tag{5.5}$$

*with a constant $0 < \rho = \rho(d, \alpha, \Gamma) < 1$.*

PROOF. Let $X \in Q_r(Y)$; we have

$$u(X) = \sum_{Z \in \partial_p Q_{2r}(Y)} u(Z)\omega^X(Z) = \sum_{Z \in \partial_p Q_{2r}(Y) - \Delta_{2r}(Y)} u(Z)\omega^X(Z).$$

Hence

$$\begin{aligned}
u(X) &\leq \omega^X[\partial_p Q_{2r}(Y) - \Delta_{2r}(Y)]M_{2r} \\
&= (1 - \omega^X[\Delta_{2r}(Y)])M_{2r} \\
&\leq (1 - \theta)M_{2r} = \rho M_{2r}.
\end{aligned}$$
$\square$

## 5.2. The Carleson principle for caloric functions vanishing on the boundary.

Let the notation be as above. Let $Y = (y, s) \in \partial_l Q$ and $c \leq r \leq R_0/2$. Assume that $u$ is a nonnegative caloric function in $Q \cap (B_{3r}(y) \times \{s - 9r^2 \leq k \leq s + 9r^2\})$ and $u = 0$ on $\partial_l Q \cap (B_{2r}(y) \times \{s - 4r^2 \leq k \leq s + 4r^2\})$. Then

$$(5.6) \qquad u(X) \leq Cu(\overline{Y}_r), \qquad X \in Q_r(Y),$$

where $C = C(d, \alpha, \Gamma) > 0$. To prove (5.6) we first observe that the local Harnack principle allows us to assume that the parabolic distance of $X$ from $\partial_p Q_{2r}(Y)$ is sufficiently large. We shall denote by $\delta(X)$ [$X \in C_{2r}(Y)$] this distance and suppose that $\delta(X) = \text{Dist}(X, \partial_p Q_{2r}(Y)) \geq C$. Geometric considerations in combination with the parabolic Harnack principle imply that

$$(5.7) \qquad \delta^\gamma(X) u(X) \leq Cr^\gamma u(\overline{Y}_r), \qquad X \in Q_{2r}(Y),$$

where $\gamma$ and $C > 0$ are positive constants depending on $d$, $\alpha$, $\Gamma$. Let $\tilde{\delta}(X) = \text{Dist}(X, \partial_p C_{2r}(Y))$. Let $0 < \varepsilon_0 < 10^{-2}$ small enough so that

$$(5.8) \qquad \theta_0 = \frac{\rho}{(1 - 4\varepsilon_0)^\gamma} < 1,$$

where $\rho$ is the constant given by Corollary 1. We shall distinguish two cases. First assume that $\delta = \delta(X) \leq \varepsilon_0 \tilde{\delta}(X)$. In this case we have $\delta = \delta(X) = \delta(x, t) = \text{dist}(x, \partial\Omega) = |x - x_0|$ for some $x_0 \in \partial\Omega$. By Corollary 1 applied to $u$ in $Q_{2\delta}(X_0) = C_{2\delta}(x_0, t) \cap Q$, we have

$$u(X) \leq \sup_{Q_{(3/2)\delta}(X_0)} u \leq \rho \sup_{Q_{3\delta}(X_0)} u.$$

But

$$\tilde{\delta}(X) \leq \tilde{\delta}(Z) + 4\delta \leq \tilde{\delta}(Z) + 4\varepsilon_0 \tilde{\delta}(X)$$

where $Z = (z, \tau) \in Q_{3\delta}(X_0)$ is such that $u(Z) = \sup_{Q_{3\delta}(X_0)} u$. Therefore

$$\tilde{\delta}(X) \leq (1 - 4\varepsilon_0)^{-1} \tilde{\delta}(Z)$$

and this gives

$$(5.9) \qquad \begin{aligned} \tilde{\delta}(X)^\gamma u(X) &\leq (1 - 4\varepsilon_0)^{-\gamma} \rho \tilde{\delta}(Z)^\gamma u(Z) \\ &\leq \theta_0 \sup_{X \in Q_{2r}(Y)} \tilde{\delta}(X)^\gamma u(X). \end{aligned}$$

It remains to examine the case where $\delta = \delta(X) > \varepsilon_0 \tilde{\delta}(X)$. In this case we have

$$(5.10) \qquad \tilde{\delta}(X)^\gamma u(X) \leq \varepsilon_0^{-\gamma} \delta(X)^\gamma u(X) \leq \varepsilon_0^{-\gamma} \sup_{Q_{2r}(Y)} \delta(X)^\gamma u(X).$$

Putting together (5.9) and (5.10) we deduce that

$$\sup_{Q_{2r}(Y)} \tilde{\delta}(X)^\gamma u(X) \leq \max\left(\theta_0 \sup_{Q_{2r}(Y)} \tilde{\delta}^\gamma u, \varepsilon_0^{-\gamma} \sup_{Q_{2r}(Y)} \delta^\gamma u\right).$$

Using (5.8), the fact that $\tilde{\delta}(X) \approx r$, $X \in Q_r(Y)$ and (5.7), we deduce the estimate (5.6).

5.3. *The boundary Harnack principle and proof of Theorem 4.* For $Y = (y, s) \in \partial\Omega \times \mathbf{Z}$; $c \leq r \leq R \leq R_0/2$, we denote

$$\Omega_{R,r}(y) = B_R(y) \cap \{x \in \Omega, \operatorname{dist}(x, \partial\Omega) < r\},$$
$$D_{R,r}(Y) = D_{R,r} = \Omega_{R,r}(y) \times \{s - R^2 \leq k \leq s\},$$
$$\Lambda_{R,r}(Y) = \Lambda_{R,r} = \partial_p D_{R,r} \cap \{x \in \Omega, 0 < \operatorname{dist}(x, \partial\Omega) < r\},$$
$$S_{R,r}(Y) = S_{R,r} = \partial_p D_{R,r} \cap \{x \in \Omega, \operatorname{dist}(x, \partial\Omega) \geq r\}.$$

LEMMA 4. *Let $Y = (y, s) \in \partial\Omega \times \mathbf{Z}$. Then we have*

(5.11) $\quad \mathbf{P}_X[\dot{S}_{\dot{\tau}_{D_{Kr,r}}} \in S_{Kr,r}] \geq \mathbf{P}_X[\dot{S}_{\dot{\tau}_{D_{Kr,r}}} \in \Lambda_{Kr,r}], \qquad X \in Q_r(Y),$

*provided that $K \geq K_0$ is large enough.*

PROOF OF THEOREM 4. Theorem 4 is an immediate consequence of the Carleson principle and Lemma 4. Indeed, we may always assume that $v(\overline{Y}_{Kr}) = u(\underline{Y}_{Kr}) = 1$. By Carleson, $v \leq c_0 = c_0(d, \alpha, \Gamma)$ in $Q_{Kr}(Y)$ (which contains $D_{Kr,r}$). The constant $c_0$ can be chosen so that, by Harnack, $u \geq 1/c_0$ on $S_{Kr,r}$. Let $u_0 = c_0 u$ and $v_0 = \frac{v}{c_0} - u_0$. Let $X \in Q_r(Y)$. By Lemma 4 we have

$$v_0(X) \leq \omega^X(\Lambda_{Kr,r}(Y)) \leq \omega^X(S_{Kr,r}(Y)) \leq u_0(X)$$

and then

$$\sup_{Q_r(Y)} \frac{v}{u} = c_0^2 \sup_{Q_r(Y)} \left(\frac{v_0}{u_0} + 1\right) \leq 2c_0^2. \qquad \square$$

PROOF OF LEMMA 4. To prove the estimate (5.11), it suffices to show that if $u, v : \overline{D_{Kr,r}} \to \mathbf{R}$ satisfy

(5.12)
$$\begin{aligned} \mathcal{L}u &= 0, & u &\geq 0 & &\text{in } D_{Kr,r}; & u &\geq 1 & &\text{on } S_{Kr,r}, \\ \mathcal{L}v &= 0, & v &\leq 1 & &\text{in } D_{Kr,r}; & v &\leq 0 & &\text{on } \partial_p D_{Kr,r} \setminus \Lambda_{Kr,r}, \end{aligned}$$

then we have

(5.13) $\qquad\qquad\qquad v \leq u \qquad \text{in } Q_r = Q_r(Y)$

provided that $K \geq K_0$ is large enough.

The first step is to prove that under (5.12), $u$ verifies the lower estimate

$$(5.14) \quad u(X) \geq 2\delta \left( \frac{\text{dist}(x, \partial \Omega)}{r} \right)^{\gamma}, \quad X = (x, t) \in Q_r(Y),$$

for appropriate constants $\delta, \gamma > 0$. Let $\tilde{y} = (R_0 - 5r) \frac{y - y_0}{|y - y_0|}$. We assume that $K \geq 10$. We define $\tilde{u} : \overline{Q_{6r}} \to \mathbf{R}$ by

$$\mathcal{L}\tilde{u} = 0 \quad \text{in } Q_{6r},$$
$$\tilde{u} = \min(u, 1) \quad \text{on } \partial_p Q_{6r} \cap \overline{D_{Kr,r}},$$
$$\tilde{u} = 1 \quad \text{on } \partial_p Q_{6r} \setminus \overline{D_{Kr,r}}.$$

Since $u \geq 0$, we have, by the maximum principle, $0 \leq \tilde{u} \leq 1$ in $Q_{6r}$, and since $u \geq 1$ on $S_{Kr,r}$ we have $u \geq \tilde{u}$ on $Q_r$. Let $\tilde{z} \in B_{2r}(\tilde{y})$ satisfying $|\tilde{z} - (R_0 - 4r) \frac{y - y_0}{|y - y_0|}| \leq 1$. Let $\tilde{Z} = (\tilde{z}, s - 4[r]^2)$. Let $w$ be defined by $w(X) = 1 - \tilde{u}(X)$, $X \in \overline{U}$, $U = Q_{6r} \cap (B_{2r}(\tilde{y}) \times \{s - 4r^2 \leq k \leq s\})$. $w$ vanishes on $\partial_l U \setminus \partial_l(B_{2r}(\tilde{y}) \times \{s - 4r^2 \leq k \leq s\})$. On the other hand, by the same argument as in Lemma 3 we see that $\omega_U^{\tilde{Z}}(\partial_l U \setminus \partial_l(B_{2r}(\tilde{y}) \times \{s - 4r^2 \leq k \leq s\})) \geq c > 0$. It follows then that $w(\tilde{Z}) \leq \theta \sup_U w$, where $0 < \theta < 1$. This means that $1 - \tilde{u}(\tilde{Z}) \leq \theta < 1$ and therefore $\tilde{u}(\tilde{Z}) \geq 1 - \theta > 0$. Using Harnack and the fact that $u \geq \tilde{u}$ on $Q_r$ we then deduce (5.14). It follows from (5.14) that

$$u(X) \geq 2\delta K^{-\gamma}, \quad X \in \overline{Q_r \setminus D_{r,r/K}},$$

where we assume $r \geq Kc$ and $K$ sufficiently large. Observe that we have in particular $u(X) \geq 2\delta K^{-\gamma}$, $X \in S_{r,r/K}$.

The second step is to prove that there exists $N > 0$ such that

$$(5.15) \quad v(X) \leq \exp(-NK), \quad X \in D_{r,r}.$$

Let $j = 1, 2, \ldots$ such that $2j + 1 \leq K$ and $X_j$ such that

$$X_j \in \partial_p D_{(2j-1)r, r}, \quad \sup_{D_{(2j-1)r, r}} v = v(X_j) = v(x_j, s_j).$$

Let $\tilde{U} = (B_{2r}(x_j) \times \{s_j - 8r^2 \leq k \leq s_j\}) \cap D_{Kr,r}$. We have $v \leq 0$ on $\partial_p U \setminus \partial_p(B_{2r}(x_j) \times \{s_j - 8r^2 \leq k \leq s_j\})$ and using the fact that

$$\omega_{\tilde{U}}^{X_j}(\partial_p \tilde{U} - \partial_p(B_{2r}(x_j) \times \{s_j - 8r^2 \leq k \leq s_j\})) \geq c > 0$$

we deduce that

$$v(X_j) \leq \sup_{\tilde{U} \cap (B_r(x_j) \times \{s_j - 4r^2 \leq k \leq s_j\})} v \leq \theta \sup_{\tilde{U}} v,$$

where $0 < \theta < 1$. Hence

$$\sup_{D_{(2j-1)r, r}} v \leq \theta \sup_{\tilde{U}} v \leq \rho \sup_{D_{(2j+1)r, r}} v,$$

where $0 < \rho < 1$. Iterating this estimate we obtain

$$\sup_{D_{r,r}} v \leq \rho^k \sup_{D_{(2k+1)r,r}} v \leq e^{-NK},$$

where $2k + 1 \leq K \leq 2k + 3$. Thus (5.15) is proved. It follows in particular that $v \leq \delta K^{-\gamma}$ in $D_{r,r}$ provided that $K$ is large enough.

From the previous considerations it follows that $u_1 = \frac{K^\gamma}{2\delta} u \geq 0$ in $D_{r,r/K}$ and $u_1 \geq 1$ on $\overline{D_{r,r} \setminus D_{r,r/K}}$ (that contains $S_{r,r/K}$) and $v_1 = \frac{K^\gamma}{2\delta}(2v - u) \leq \frac{K^\gamma}{\delta} v \leq 1$ in $D_{r,r}$ (that contains $D_{r,r/K}$) with $v_1 \leq 0$ on $S_{r,r/K}$. In particular, we have

$$u_1 - v_1 = \frac{K^\gamma}{\delta}(u - v) \geq 0, \text{ on } \overline{D_{r,r} \setminus D_{r,r/K}}.$$

On the other hand, $u_1$, $v_1$ satisfy the same assumptions as $u$, $v$ with $r$ replaced by $r/K$. We can iterate and define $u_j$, $v_j$ such that

$$u_j - v_j = \left(\frac{K^\gamma}{\delta}\right)^j (u - v) \geq 0 \qquad \text{on } \overline{D_{r/K^j, r/K^j} \setminus D_{r/K^j, r/K^{j+1}}},$$

$j = 1, 2, \ldots$, and consequently

$$u - v \geq 0 \qquad \text{on } S(Y) = \bigcup_{j \geq 0} \overline{D_{r/K^j, r/K^j} \setminus D_{r/K^j, r/K^{j+1}}}.$$

Let now $X_0 = (x_0, t_0) \in Q_r \subset D_{r,r}(Y)$ and let $\tilde{X}_0 = (\tilde{x}_0, t_0)$ where $\tilde{x}_0 \in \partial\Omega$ satisfies $\text{dist}(x, \partial\Omega) = |x_0 - \tilde{x}_0|$. Then $D_{Kr,r}(\tilde{X}_0) \subset D_{(K+2)r,r}(Y)$ and $S_{Kr,r}(\tilde{X}_0) \subset S_{(K+2)r,r}(Y)$. Replacing $K$ with $K + 2$ in the previous considerations, we deduce that $u \geq v$ on $S(\tilde{X}_0)$ that contains $X_0$. This completes the proof of Lemma 4. □

5.4. *The boundary backward Harnack principle.* Let the notation be the same as in Theorem 5. Let $B_r(y_0) = B_r$. First we observe that the Carleson principle in combination with the parabolic Harnack principle give

(5.16) $\qquad u\left(x, \left[\frac{r^2}{8}\right]\right) \leq C \min_{B_{r/2} \times \{r^2/4 \leq t \leq 8r^2\}} u, \qquad x \in B_r - B_{r/2}.$

On the other hand, by Harnack

(5.17) $\qquad \max_{B_{r/2} \times \{t=[r^2/8]\}} u \leq C \min_{B_{r/2} \times \{r^2/4 \leq t \leq 8r^2\}} u.$

Since $u \equiv 0$ on $\partial B_r \times \mathbf{N}$, by (5.16), (5.17) and the maximum principle we get

$$\max_{B_r \times \{[r^2/8] \leq t \leq 8r^2\}} u \leq C \min_{B_{r/2} \times \{r^2/4 \leq t \leq 8r^2\}} u.$$

In particular, we have

(5.18) $$\max_{B_{r/2} \times \{[r^2/8] \leq t \leq 8r^2\}} u \leq C \min_{B_{r/2} \times \{r^2/4 \leq t \leq 8r^2\}} u.$$

Let now $v$ be another nonnegative caloric function in $B_r \times \mathbf{N}$ such that $v$ vanishes on the lateral boundary $\partial B_r \times \mathbf{N}$. Then, there exists $C > 0$ such that

(5.19) $$u(x_0, [r]^2)v(x,t) \leq Cv(x_0, 4[r]^2)u(x,t),$$
$$x_0 \in B_{r/2}, (x,t) \in B_r \times \{r^2 \leq t \leq 3r^2\}.$$

To prove (5.19) we first use a covering argument to cover $\partial_l(B_r \times \{r^2 \leq t \leq 3r^2\})$ with cylinders $C_j = B_{\varepsilon r}(y_j) \times \{s_j - \varepsilon^2 r^2 \leq t \leq s_j\}$, $j = 1, \ldots, N$, where $\varepsilon > 0$ is chosen sufficiently small and where $Y^{(j)} = (y_j, s_j) \in \partial_l(B_r \times \{r^2 \leq t \leq 3r^2\})$, and apply in each of these $C_j$ the boundary Harnack principle to get

(5.20) $$u(\underline{Y}^{(j)}_{r/8})v(x,t) \leq Cv(\overline{Y}^{(j)}_{r/8})u(x,t), \qquad (x,t) \in C_j,$$

where $\underline{Y}^{(j)}_r$, $\overline{Y}^{(j)}_r$, $r > 0$, are defined as in Section 2. By Harnack, we have

(5.21) $$v(\overline{Y}^{(j)}_{r/8}) \leq Cv(x_0, 4[r]^2),$$

(5.22) $$u(\underline{Y}^{(j)}_{r/8}) \geq cu(x_0, [r^2/4]).$$

Using (5.18) we deduce from (5.22) that

(5.23) $$u(\underline{Y}^{(j)}_{r/8}) \geq cu(x_0, [r]^2).$$

Putting together (5.20), (5.21) and (5.23) we deduce that

(5.24) $$u(x_0, [r]^2)v(x,t) \leq Cv(x_0, 4[r]^2)u(x,t),$$
$$x_0 \in B_{r/2}, (x,t) \in C_j, j = 1, \ldots, N.$$

On the other hand, we have, by Harnack,

(5.25) $$v(x,t) \leq Cv(x_0, 4[r]^2), \qquad r^2 \leq t \leq 3r^2, \qquad \text{dist}(x, \partial B_r) \geq \delta r,$$

where $0 < \delta < 1/2$. Again, by Harnack combined with (5.18)

(5.26) $$u(x,t) \geq cu(x_0, [r]^2), \qquad r^2 \leq t \leq 3r^2, \qquad \text{dist}(x, \partial B_r) \geq \delta r$$

and (5.19) follows from (5.24)–(5.26) with an appropriate choice of $\delta > 0$. We are now able to get the estimate (2.3). We shall use (5.18), (5.19) and a time-shifting argument. Let $u$ be as in (2.3) and let $v(x,t) = u(x, t + 2[r]^2)$. By (5.19) we have

(5.27) $$u(x_0, [r]^2)u(x, t + 2[r]^2) \leq Cu(x_0, 6[r]^2)u(x,t),$$
$$(x,t) \in B_r \times \{r^2 \leq t \leq 3r^2\},$$

with $x_0 \in B_{r/2}$ fixed. Equation (2.3) follows from (5.27) and the estimate

$$u(x_0, 6[r]^2) \leq Cu(x_0, [r]^2),$$

which is an immediate consequence of (5.18). This completes the proof of Theorem 5.

## REFERENCES

[1] ARONSON, D. G. (1968). Non-negative solutions of linear parabolic equations. *Ann. Sci. Norm. Super. Pisa* **22** 607–694. MR0435594
[2] AUSCHER, P. and COULHON, T. (1999). Gaussian lower bounds for random walks from elliptic regularity. *Ann. Inst. H. Poincaré Probab. Statist.* **35** 605–630. MR1705682
[3] BASS, R. F. and BURDZY, K. (1994). The boundary Harnack principle for nondivergence form elliptic operators. *J. London Math. Soc.* **50** 157–169. MR1277760
[4] BAUMAN, P. E. (1984). Positive solutions of elliptic equations in nondivergence form and their adjoints. *Ark. Mat.* **22** 536–565. MR0765409
[5] DELMOTTE, T. (1999). Parabolic Harnack inequality and estimates of Markov chains on graphs. *Rev. Mat. Iberoamericana* **15** 181–232. MR1681641
[6] ESCAURIAZA, L. (2000). Bounds for the fundamental solution of elliptic and parabolic equations in nondivergence form. *Comm. Partial Differential Equations* **25** 821–845. MR1759794
[7] FABES, E. B., GAROFALO, N. and SALSA, S. (1986). A backward Harnack inequality and Fatou theorem for nonnegative solutions of parabolic equations. *Illinois J. Math.* **30** 536–565. MR0857210
[8] FABES, E. B. and SAFONOV, M. V. (1997). Behavior near the boundary of positive solutions of second order parabolic equations. *J. Fourier Anal. Appl.* **3** 871–882. MR1600211
[9] FABES, E. B., SAFONOV, M. V. and YUAN, Y. (1999). Behavior near the boundary of positive solutions of second order parabolic equations. II. *Trans. Amer. Math. Soc.* **351** 4947–4961. MR1665328
[10] FABES, E. B. and STROOCK, D. W. (1986). A new proof of Moser's parabolic Harnack inequality via the old idea of Nash. *Arch. Rational Mech. Anal.* **96** 327–338. MR0855753
[11] GRIGOR'YAN, A. (1991). The heat equation on noncompact Riemannian manifolds. *Mat. Sb.* **182** 55–87. [Translation in *Russian Acad. Sci. Sb. Math.* **72** (1992) 47–77.] MR1098839
[12] KOZLOV, S. M. (1985). The method of averaging and random walks in inhomogeneous environments. *Russian Math. Surveys* **40** 73–145.
[13] KUO, H. J. and TRUDINGER, N. S. (1998). Evolving monotone difference operators on general space–time meshes. *Duke Math. J.* **91** 587–607. MR1604175
[14] LAWLER, G. F. (1992). Estimates for differences and Harnack inequality for difference operators coming from random walks with symmetric, spatially inhomogeneous, increments. *Proc. London Math. Soc.* **63** 552–568. MR1127149
[15] SAFONOV, M. V. and YUAN, Y. (1999). Doubling properties for second order operators. *Ann. of Math.* (*2*) **150** 313–327. MR1715327
[16] SALOFF-COSTE, L. (1995). Parabolic Harnack inequality for divergence-form second-order differential operators. Potential theory and degenerate partial differential operators (Parma). *Potential Anal.* **4** 429–467. MR1354894

[17] SALOFF-COSTE, L. and HEBISCH, W. (1993). Gaussian estimates for Markov chains and random walks on groups. *Ann. Probab.* **21** 673–709. MR1217561
[18] VAROPOULOS, N. TH. (2000). Potential theory in conical domains. II. *Math. Proc. Cambridge Philos. Soc.* **129** 301–319. MR1765917
[19] VAROPOULOS, N. TH., SALOFF-COSTE, L. and COULHON, T. (1992). *Analysis and Geometry on Groups.* Cambridge Univ. Press. MR1218884

INSTITUT MATHÉMATIQUE DE JUSSIEU
175 RUE DU CHEVALERET
75013 PARIS
FRANCE
E-MAIL: sam@math.jussieu.fr